\documentclass[oneside, 12pt]{amsart} 

\usepackage[utf8]{inputenc}
\usepackage[margin=1in]{geometry}
\usepackage{amsthm}
\usepackage{enumitem}  
\usepackage[OT2,T1]{fontenc}
\usepackage{comment}
\usepackage{amscd,amssymb, enumitem, amsmath,mathrsfs, amsfonts, amsaddr}
\usepackage{url}
\usepackage{tikz}
\usepackage{fullpage}
\usepackage{microtype}
\usetikzlibrary{decorations.pathreplacing}
\usepackage[backref=page]{hyperref}
\usepackage{hyperref}
\usepackage[utf8]{inputenc}

\DeclareSymbolFont{cyrletters}{OT2}{wncyr}{m}{n}
\DeclareMathSymbol{\Sha}{\mathalpha}{cyrletters}{"58}

\usepackage[utf8]{inputenc}

\newtheorem{theorem}{Theorem}[section]

\newtheorem{proposition}[theorem]{Proposition}
\newtheorem*{proposition*}{Proposition}

\newtheorem*{questiona*}{Question A}
\newtheorem*{questionb*}{Question B}
\newtheorem*{theorem*}{Theorem}

\theoremstyle{definition}

\newtheorem{example}[theorem]{Example}

\newtheorem{remark}[theorem]{Remark}
\newtheorem{emptyremark}[theorem]{}
\newtheorem*{acknowledgement}{Acknowledgements}

\theoremstyle{remark}

\title{Tamagawa numbers of elliptic curves with torsion points}
\author{Mentzelos Melistas}
\address{ Steklov Mathematical Institute of Russian Academy of Sciences, Moscow, Russia,\\ ~~~~~ email: mentzmel@gmail.com}
\date{\today}

\begin{document}

\maketitle

\begin{abstract}
Let $K$ be a global field and let $E/K$ be an elliptic curve with a $K$-rational point of prime order $p$. In this paper we are interested in how often the (global) Tamagawa number $c(E/K)$ of $E/K$ is divisible by $p$. This is a natural question to consider in view of the fact that the fraction $c(E/K)/ |E(K)_{\text{tors}}|$ appears in the second part of the Birch and Swinnerton-Dyer Conjecture. We focus on elliptic curves defined over global fields, but we also prove a result for higher dimensional abelian varieties defined over $\mathbb{Q}$.
\end{abstract}

\section{Introduction}

Let $K$ be a global field and let $E/K$ be an elliptic curve. Let $v$ be a non-archimedean valuation of $K$. We denote by $K_{v}$ the completion of $K$ with respect to the valuation $v$, by $\mathcal{O}_{K_{v}}$ the valuation ring of $K_{v}$, and by $k_{v}$ the residue field of $K_{v}$. The set $E_0(K_v)$ that consists of points with nonsingular reduction is a finite index subgroup of $E(K_v)$. The index $c_{v}(E/K)=[E(K_v):E_0(K_v)]$ is called the Tamagawa number of $E/K$ at $v$. Alternatively, the number $c_{v}(E/K)$ can be defined as $|\mathcal{E}_{k_v}(k_v)/\mathcal{E}^0_{k_v}(k_v)|$, where $\mathcal{E}_{k_v}/k_v$ is the special fiber of the N\'eron model of $E/K$ and $\mathcal{E}^0_{k_v}/k_v$ is the connected component of the identity of $\mathcal{E}_{k_v}/k_v$. The two above definitions of $c_{v}(E/K)$ agree by  \cite[Corollary IV.9.2]{silverman2}. We define the (global) Tamagawa number of $E/K$ as $c(E/K):=\prod_{v} c_{v}(E/K)$, where the product is taken over all the non-archimedean valuations of $K$. If $\mathfrak{p}$ be a prime of the ring of integers $\mathcal{O}_K$ of $K$ corresponding to a non-archimedean valuation $v$, then we will also denote $c_v(E/K)$ by $c_{\mathfrak{p}}(E/K)$.

Let $K$ be a global field and let $E/K$ be an elliptic curve (or more generally an abelian variety). In this paper, we are interested in the effect of the torsion subgroup of $E/K$ to the Tamagawa number of $E/K$. The importance of such results stems from the fact that the fraction $c(E/K)/ |E(K)_{\text{tors}}|$ appears in the Birch and Swinnerton-Dyer Conjecture (see \cite[Appendix C.16]{aec} or \cite[Conjecture F.4.1.6]{hindrysilverman}). 

The first to take up the study of Tamagawa numbers of abelian varieties with torsion points was Lorenzini in \cite{lor}. Krumm in \cite{Krummthesis} investigated the interplay between torsion points and Tamagawa numbers for elliptic curves over number fields of low degree and formulated a conjecture concerning elliptic curves with a point of order $13$ over quadratic number fields. Krumm's conjecture was later proved by Najman in \cite{najmantamawanumber}. Very recently, Trbović has studied in \cite{trbovic2021tamagawa} Tamagawa numbers of elliptic curves with an isogeny.

When $K$ is a quadratic number field, Lorenzini (see \cite[Corollary 3.4]{lor}) has proved that there are no elliptic curves $E/K$ with a $K$-rational point of order $11$ and such that $11 \nmid c(E/K)$. Krumm (see \cite[Proposition 4.2]{Krummthesis}) proved the same result when $K$ is a cubic number field. Example \ref{example11torsionoverdegree4} in the next section shows that there exists a quartic number field $K$ and an elliptic curve $E/K$ with a $K$-rational point of order $11$ such that $c(E/K)=1$. Since the modular curve $X_1(11)$ has genus $1$, for a given number field $K$, there may exist infinitely many elliptic curves $E/K$ with a $K$-rational point of order $11$. Theorem \ref{theorem11torsionpoint}, which is proved in Section \ref{proofsofresults}, provides a generalization of Lorenzini and Krumm's results. 

\begin{theorem}\label{theorem11torsionpoint}
For every number field $K/\mathbb{Q}$ there exists a constant $n_{K,11}$ such that the following holds: For every elliptic curve $E/K$ with a $K$-rational point of order $11$ we have that $11$ divides $c(E/K)$ with at most $n_{K,11}$ exceptions.
\end{theorem}

Let now $K$ be equal to $\mathbb{Q}$ or a quadratic number field. It follows from work of Lorenzini (see \cite[Proposition 1.1, Proposition 2.10, and Corollary 3.4]{lor}) that if $E/K$ is an elliptic curve with a $K$-rational point of prime order $p \geq 7$, then $p$ divides $c(E/K)$ with only finitely many exceptions (the number of exceptions depends on $K$). The following theorem, which follows from Theorem \ref{theorem11torsionpoint} and is proved in Section \ref{proofsofresults}, provides a generalization of the above statement to number fields of arbitrary degree.

\begin{theorem}\label{firsttheorem}
For every number field $K/\mathbb{Q}$ there exists a constant $n_K$ such that the following holds: For every prime $p \geq 7$ and every elliptic curve $E/K$ with a $K$-rational point of order $p$ we have that $p$ divides $c(E/K)$ with at most $n_K$ exceptions.
\end{theorem}

It follows from \cite[Lemma 2.26]{lor} and \cite[Corollary 5.4]{barriosroy} that for $p=2$ or $3$ there exist infinitely many elliptic curves $E/\mathbb{Q}$ with a $\mathbb{Q}$-rational point of order $p$ such that $c(E/\mathbb{Q})=1$. Also, as we explain in Remark \ref{remarkp235}, it seems likely that there exists a number field $K/\mathbb{Q}$ and an infinite number of elliptic curves $E/K$ with a $K$-rational point of order $5$ and such that $5 \nmid c(E/\mathbb{Q})$. Therefore, the assumption that $p \geq 7$ is necessary in Theorem \ref{firsttheorem}. Moreover, Examples \ref{example1Krumm} and \ref{example2Krumm} below, which are due to Krumm, show that there exist elliptic curves $E/K$ defined over cubic and quartic number fields $K/\mathbb{Q}$ with a $K$-rational point of order greater than $11$ and such that $c(E/K)=1$.

One may wonder whether Theorem \ref{theorem11torsionpoint} can be generalized to abelian varieties of higher dimension by requiring that the constant also depends on their dimension. As Part $(i)$ of the following theorem shows, such a generalization is not possible even for $K=\mathbb{Q}$.

\begin{theorem}\label{weilrestrictiontamagawa}
Let $p \geq 5$ be a prime. 
 \begin{enumerate}[topsep=2pt,label=(\roman*)]
\itemsep0em 
  \item There exist infinitely many abelian varieties $A/\mathbb{Q}$ of dimension at most $\frac{p^2-1}{2}$ with a $\mathbb{Q}$-rational point of order $p$ and such that $p \nmid c(A/\mathbb{Q})$.
  \item If $p \equiv 1 \; (\text{mod} \; 3)$, then there exists an abelian variety $A/\mathbb{Q}$ of dimension $\frac{p-1}{3}$ with a $\mathbb{Q}$-rational point of order $p$ and such that $p \nmid c(A/\mathbb{Q})$.
  \item If $p \equiv 2 \; (\text{mod} \; 3)$ and $p \equiv 1 \; (\text{mod} \; 4)$, then there exists an abelian variety $A/\mathbb{Q}$ of dimension $\frac{p-1}{2}$ with a $\mathbb{Q}$-rational point of order $p$ and such that $p \nmid c(A/\mathbb{Q})$.
  \end{enumerate}
\end{theorem}

\begin{remark}
Let $p \geq 5$ be a prime and let $f_p $ be the minimum out of all $d>0$ that satisfy the following statement: there exist only finitely many abelian varieties $A/\mathbb{Q}$ of dimension $d$ with a $\mathbb{Q}$-rational point of order $p$ and such that $p \nmid c(A/\mathbb{Q})$. It follows from \cite[Proposition 1.1]{lor} combined a celebrated theorem of Mazur (see \cite[Theorem (8)]{maz}) on the classification of all the possible rational torsion subgroups of rational elliptic curves that $1 \leq f_p$. Theorem \ref{weilrestrictiontamagawa} shows that $f_p \leq \frac{p^2-1}{2}$.
\end{remark}

On the other hand, keeping the same notation as above, if we let $p$ depend on the dimension of the abelian variety as well as the degree of the base field, then Tamagawa number divisibility by the prime $p$ can be achieved. More precisely, it follows from \cite[Proposition 3.1]{lor} that given a number field $K$ and an integer $d>0$, then there exists a constant $\gamma_{K,d}$ such that if $A/K$ is an abelian variety with a $K$-rational point of order $p \geq \gamma_{K,d}$, then $p$ divides $c(A/K)$.

One can also study Tamagawa numbers of elliptic curves defined over function fields. Here we will focus on the case where the elliptic curve is defined over a function field of characteristic $p$ and has a point of order $p$. More specifically, in Section \ref{sectionfunctionfields} we prove the following theorem.

\begin{theorem}
Let $p \geq 5$ be a prime and let $q$ be a power of $p$.
\begin{enumerate}
    \item Let $E/\mathbb{F}_q(t)$ be a non-isotrivial elliptic curve with an $\mathbb{F}_q(t)$-rational point of order $p$. Then $p$ divides $c(E/\mathbb{F}_q(t))$.
    \item Let $K=\mathbb{F}_q(\mathcal{C})$ be the function field of a smooth, projective, geometrically irreducible curve $\mathcal{C}/\mathbb{F}_q$ and let $E/K$ be a non-isotrivial elliptic curve with a $K$-rational point of order $p$. Then there exists a finite extension $K'/K$ such that $p$ divides $c(E_{K'}/K')$, where $E_{K'}/K'$ is the base change of $E/K$ to $K'$.
\end{enumerate}
\end{theorem}

\begin{acknowledgement}
This project was initiated while the author was a graduate student at the University of Georgia and a part of this work is contained in the author’s doctoral dissertation. The author would like to thank Dino Lorenzini for many useful suggestions during the preparation of this work and Pete Clark for providing an argument that improved a previous version of Theorem \ref{weilrestrictiontamagawa}. I would like to thank the anonymous referee for many insightful comments and many useful suggestions. This work was completed at the Steklov International Mathematical Center and supported by the Ministry of Science and Higher education of the Russian Federation (Agreement no. 075-15-2019-1614).
\end{acknowledgement}

\section{Proofs of Theorems \ref{theorem11torsionpoint}, \ref{firsttheorem}, and \ref{weilrestrictiontamagawa}}\label{proofsofresults}
 
\begin{emptyremark}\label{tatealgorithm}
Tate, in \cite{tatealgorithm}, has produced an algorithm that computes the Tamagawa number of an elliptic curve defined over a complete discrete valuation ring. We recall a part of the algorithm that we will use. We refer the reader to \cite[Section IV.9]{silverman2} (or \cite{tatealgorithm}) for more details. Let $R$ be a complete discrete valuation ring with (normalized) valuation $v$, fraction field $K$, and perfect residue field $k$. Let $E/K$ be an elliptic curve given by a Weierstrass equation $$y^2+a_1xy+a_3y=x^3+a_2x^2+a_4x+a_6,$$ with $a_i \in R$ for $i = 1,2,3,4,6,$ and such that $v(c_4)=0$ and $v(\Delta) > 0$. Here $\Delta$ is the discriminant and $c_4$ is the $c_4$-invariant of the Weierstrass equation. Since $v( \Delta) >0$, we can make a change of variables so that $v(a_3), v(a_4), v(a_6) >0$. Set $b_2:=a_1^2+4a_2$. Since $v(c_4) = 0$ and $v(a_3), v(a_4), v(a_6) >0$, we obtain that $v(b_2)=0$ (see \cite[Section III.1]{aec} for the standard equation involving $c_4$ and $b_2$). Let $k'$ be the splitting field over $k$ of the polynomial $T^2+a_1T+a_2$. The curve $E/K$ has split multiplicative reduction of type I$_n$ if $v(\Delta)=n$ and $k'=k$.  In this case, the Tamagawa number $c_v(E/K)$ of $E/K$ at $v$ is equal to $n$. We will use the following observation repeatedly in this paper.

\underline{Observation:} If $v(a_2), v(a_3), v(a_4), v(a_6) >0,$ and $v(a_1)=0$, then $E/K$ has split multiplicative reduction of type I$_n$, where $v(\Delta)=n$. Moreover, in this case $n = c_v(E/K)$.

\end{emptyremark}

\begin{proof}[Proof of Theorem \ref{theorem11torsionpoint}]

 Let $K$ be a number field and let $E/K$ be an elliptic curve with a $K$-rational point of order $11$. The curve $E/K$ can be given by an equation of the form
\begin{align}\label{11torsiongeneralform}
   E(r,s): y^2+(1-c)xy-by=x^3-bx^2,
\end{align}
where $$b=rs(r-1), \quad c=s(r-1),$$
and $r,s$ satisfy $$F_{11}(r,s)\; : \; r^2 - rs^3 + 3rs^2 - 4rs + s=0.$$
The above equation is the raw form of the affine modular curve $Y_1(11)$ (see \cite{sutherland} and \cite{sutherlandrawforms}). \stepcounter{theorem}

Using SAGE \cite{sagemath} we find that the discriminant of $E(r,s)$ is $$\Delta_{E(r,s)}=r^3s^4(r - 1)^5(r^2s^3 - 8r^2s^2 - 2rs^3 + 16r^2s + 5rs^2 + s^3 - 20rs + 3s^2 + 3s + 1).$$
Let $$f(r,s)=r^2s^3 - 8r^2s^2 - 2rs^3 + 16r^2s + 5rs^2 + s^3 - 20rs + 3s^2 + 3s + 1.$$

If $\mathfrak{p}$ is any prime of $\mathcal{O}_K$, then we denote by $v_{\mathfrak{p}}$ the valuation of $\mathcal{O}_K$ associated to $\mathfrak{p}$.

\begin{proposition}\label{11torsionvaluationofr}
Let $\mathfrak{p}$ be a prime of $\mathcal{O}_K$. If $v_{\mathfrak{p}}(r) \neq 0$, then $E/K$ has split multiplicative reduction modulo $\mathfrak{p}$ and $11 \mid c_{\mathfrak{p}}(E/K)$.
\end{proposition}
\begin{proof}
We split the proof into cases. Suppose first that $v_{\mathfrak{p}}(r)>0$ and $v_{\mathfrak{p}}(s)<0$. This implies that $v_{\mathfrak{p}}(r-1)=0$. By looking at the expression for $F_{11}(r,s)$ and keeping in mind that that $F_{11}(r,s)=0$, we must have that there are two terms of the same minimal valuation. Therefore, we see that $v_{\mathfrak{p}}(s)=v_{\mathfrak{p}}(rs^3)$ and, hence, $v_{\mathfrak{p}}(r)=-2v_{\mathfrak{p}}(s)$. By performing a change of variables of the form $(x,y) \rightarrow{(s^2 x, s^3 y)}$ in Equation $($\ref{11torsiongeneralform}$)$ we obtain the following Weierstrass equation
$$y^2+(\frac{1}{s}-(r-1))xy+\frac{1}{s^2}r(r-1)y=x^3-\frac{1}{s}r(r-1)x^2.$$
By looking at the valuation of each coefficient, using the observation of \ref{tatealgorithm}, we see that the above equation is a minimal Weierstrass equation and moreover that $E/K$ has split multiplicative reduction modulo $\mathfrak{p}$. The discriminant of the new equation is $${\Delta'}_{E(r,s)}=\frac{r^3s^4(r-1)^5f(r,s)}{s^{12}}=\frac{r^3(r-1)^5f(r,s)}{s^8}$$ and since $v_{\mathfrak{p}}(f(r,s))=3v_{\mathfrak{p}}(s)$, we have that
$$v_{\mathfrak{p}}({\Delta'}_{E(r,s)})=3v_{\mathfrak{p}}(r)+3v_{\mathfrak{p}}(s)-8v_{\mathfrak{p}}(s)=3v_{\mathfrak{p}}(r)-5v_{\mathfrak{p}}(s)=-6v_{\mathfrak{p}}(s)-5v_{\mathfrak{p}}(s)=-11v_{\mathfrak{p}}(s).$$

Suppose that $v_{\mathfrak{p}}(r)>0$ and $v_{\mathfrak{p}}(s) \geq 0$. This implies that $v_{\mathfrak{p}}(r-1)=0$. By looking at the valuation of each term in the expression for $F_{11}(r,s)$, we find that $2v_{\mathfrak{p}}(r)=v_{\mathfrak{p}}(s)$. By Equation $($\ref{11torsiongeneralform}$)$, using the observation of \ref{tatealgorithm}, we see that $E/K$ has split multiplicative reduction modulo $\mathfrak{p}$ with
$$v_{\mathfrak{p}}(\Delta_{E(r,s)})=3v_{\mathfrak{p}}(r)+4v_{\mathfrak{p}}(s)=3v_{\mathfrak{p}}(r)+8v_{\mathfrak{p}}(r)=11v_{\mathfrak{p}}(r)$$

Suppose now that $v_{\mathfrak{p}}(r)<0$ and $v_{\mathfrak{p}}(s)<0$. This implies that $v_{\mathfrak{p}}(r-1)=v_{\mathfrak{p}}(r)$. By looking at the expression for $F_{11}(r,s)$ we see that $v_{\mathfrak{p}}(r^2)=v_{\mathfrak{p}}(rs^3)$ and, therefore, $v_{\mathfrak{p}}(r)=3v_{\mathfrak{p}}(s)$. By performing a change of variables of the form $(x,y) \rightarrow{(s^2(r-1)^2 x, s^3(r-1)^3 y)}$ in Equation $($\ref{11torsiongeneralform}$)$ we get a new Weierstrass equation $$y^2+\Big(\frac{1}{s(r-1)}-1\Big)xy+\frac{r}{s^2(r-1)^2}y=x^3-\frac{r}{s(r-1)}x^2.$$
This equation is an integral Weierstrass equation and moreover $E/K$ has split multiplicative modulo $\mathfrak{p}$ by the observation of \ref{tatealgorithm}. The discriminant of the new equation is 
$${\Delta'}_{E(r,s)}=\frac{r^3s^4(r-1)^5f(r,s)}{s^{12}(r-1)^{12}}=\frac{r^3f(r,s)}{s^8(r-1)^7}$$ and moreover $v_{\mathfrak{p}}(f(r,s))=2v_{\mathfrak{p}}(r)+3v_{\mathfrak{p}}(s)$. Therefore, we obtain that
$$v_{\mathfrak{p}}({\Delta'}_{E(r,s)})=3v_{\mathfrak{p}}(r)-8v_{\mathfrak{p}}(s)-7v_{\mathfrak{p}}(r)+2v_{\mathfrak{p}}(r)+3v_{\mathfrak{p}}(s)=-2v_{\mathfrak{p}}(r)-5v_{\mathfrak{p}}(s)=-6v_{\mathfrak{p}}(s)-5v_{\mathfrak{p}}(s)=-11v_{\mathfrak{p}}(s).$$

Suppose now that $v_{\mathfrak{p}}(r)<0$ and $v_{\mathfrak{p}}(s) \geq 0$. This case is impossible because by looking at the expression for $F_{11}(r,s)$, we find that $v_{\mathfrak{p}}(F_{11}(r,s))=2v_{\mathfrak{p}}(r)<0$ but $F_{11}(r,s)=0$.
\end{proof}

We are now ready to complete the proof of Theorem \ref{theorem11torsionpoint}. Let  $$T= \{ (r,s) \; : \; r,s \in \mathcal{O}_K^\ast \text{ and } F_{11}(r,s)=0 \}.$$

Let $E/K$ be an elliptic curve given by parameters $r$ and $s$. We first show that $11 \mid c(E/K)$ except (possibly) for $(r,s) \in T$. Proposition \ref{11torsionvaluationofr} implies that if $v_{\mathfrak{p}}(r) \neq 0$ for some prime $\mathfrak{p}$, then $11 \mid c_{\mathfrak{p}}(E/K)$. Therefore, if $11 \nmid c(E/K)$, then $r \in \mathcal{O}_K^\ast$. Moreover, since $r \in \mathcal{O}_K^\ast$ and $F_{11}(r,s)=0$, we obtain that $s \in \mathcal{O}_K^\ast$.

Finally, since $F_{11}(r,s)=0$ defines a (geometric) genus 1 affine curve, Siegel's Theorem implies that the set $T$ is finite (see \cite[Theorem 7.3.9]{heightsindiophantinegeometry} or \cite[Remark D.9.2.2]{hindrysilverman}). This proves our theorem.
\end{proof}

\begin{example}\label{example11torsionoverdegree4}
This example shows that there exists a number field $K/\mathbb{Q}$ of degree $4$ and an elliptic curve $E/K$ with $K$-rational point of order $11$ and such that $c(E/\mathbb{Q})=1$. Let $$K:=\mathbb{Q}[a]/(a^4-a^3-3a^2+a+1))$$ and consider the elliptic curve $E/K$ given by 
$$y^2+(a^3-3a)xy+(a^2-a)y=x^3+(-a^3+2a^2+a-3)x^2+(-a^2+1)x-a^2+a+1.$$ The curve $E/K$ is the elliptic curve with LMFDB \cite{lmfdb} label 4.4.725.1-109.1-a2, has a $K$-rational point of order $11$, and $c(E/\mathbb{Q})=1$.
\end{example}

\begin{proof}[Proof of Theorem \ref{firsttheorem}]
Let $K/\mathbb{Q}$ be a number field of degree $d$. Merel's Theorem on the boundedness of torsion of elliptic curves over number fields, see \cite[Th\'eorème]{merel1996}, implies that if $E(K)$ contains a point of order prime order $p$, then $p < d^{3d^2}$. Moreover, if  $p > 11$ is a prime, then the modular curve $X_1(p)$ has genus greater or equal to $2$. Therefore, Falting's Theorem implies that for each prime $p$ with $11<p<d^{3d^2}$ there are only finitely many elliptic curves that are defined over $K$ and have a $K$-rational point of order $p$. If $E/K$ has a rational point of order $7$, then \cite[Proposition 2.10]{lor} implies that $7 \mid c(E/K)$ with only finitely many exceptions. Therefore, in order to prove Theorem \ref{firsttheorem} it is enough to consider the case $N=11$. This is exactly Theorem \ref{theorem11torsionpoint}.
\end{proof}

\begin{remark}\label{remarkp235}
As noted in \cite[Remark 2.8]{lor}, in order to produce a number field $K/\mathbb{Q}$ and an infinite number of elliptic curves $E/K$ with a $K$-rational point of order $5$ and such that $5 \nmid c(E/K)$ it is enough to find an infinite number of units $\lambda \in \mathcal{O}_K^*$ such that the order of $\lambda^2-\lambda +1$ at any prime $\mathfrak{p}$ of $\mathcal{O}_K$ is not divisible by $5$. It seems likely that this is possible.
\end{remark}

\begin{proof}[Proof of Theorem \ref{weilrestrictiontamagawa}]
Let $p \geq 5$ be a prime.

\textit{Proof of $(i)$:} The degree of the map $\pi : X_1(p) \rightarrow X_1(1) \cong \mathbb{P}^1$ coming from the $j$-invariant is $\frac{p^2-1}{2}$ (see \cite[Page 66]{diamondshurman}).  For every $n \in \mathbb{Z}$, let $E_n/L_n$ be a closed point in the fiber of $\pi$ over $n$, i.e., the curve  $E_n/L_n$ is defined over a number field $L_n/\mathbb{Q}$ of degree at most $\frac{p^2-1}{2}$ and $j(E_n)=n$. Let $A_n/\mathbb{Q}$ be the Weil restriction of $E_{n}/L_n$ to $\mathbb{Q}$ (see \cite[Section 7.6]{neronmodelsbook} for the basics of Weil restriction). Since $j(E_n)=n \in \mathbb{Z}$, the curve $E_{n}/L_n$ has everywhere potentially good reduction and, hence, $p \nmid c(E_{n}/L_{n})$ because $p \geq 5$. It follows from \cite[Proposition 3.19]{lor} that $c(A_n/\mathbb{Q})=c(E_{n}/L_n)$. Therefore, for every $n$ the abelian variety $A_n/\mathbb{Q}$ has a $\mathbb{Q}$-rational point of order $p$, has dimension at most $\frac{p^2-1}{2}$, and $p \nmid c(A_n/\mathbb{Q})$. Thus part $(i)$ is proved.

\textit{Proof of $(ii)$:} Since $p \equiv 1 \; (\text{mod} \; 3)$, Part $(a)$ of \cite[Theorem 1]{clarkcookstankewicz} implies that there exist a field extension $K/\mathbb{Q}$ of degree $\frac{p-1}{3}$ and a CM elliptic curve $E/K$ with a $K$-rational point of order $p$. Since $p \geq 5$ and $E/K$ has potentially good reduction, we find that $p \nmid c(E/K)$. Let $A/\mathbb{Q}$ be the Weil restriction of $E/K$ to $\mathbb{Q}$. It follows from \cite[Proposition 3.19]{lor} that $c(A/\mathbb{Q})=c(E/K)$. Therefore, the abelian variety $A/\mathbb{Q}$ has $\mathbb{Q}$-rational point of order $p$, has dimension $\frac{p-1}{3}$, and $p \nmid c(A/\mathbb{Q})$. This proves part $(ii)$.

\textit{Proof of $(iii)$:} The proof is similar to the proof of part $(ii)$. Since  $p \equiv 2 \; (\text{mod} \; 3)$ and $p \equiv 1 \; (\text{mod} \; 4)$, Part $(b)$ of \cite[Theorem 1]{clarkcookstankewicz} implies that there exist a field extension $K/\mathbb{Q}$ of degree $\frac{p-1}{2}$ and a CM elliptic curve $E/K$ with a $K$-rational point of order $p$.  Since $p \geq 5$ and $E/K$ has potentially good reduction, we find that $p \nmid c(E/K)$. Let $A/\mathbb{Q}$ be the Weil restriction of $E/K$ to $\mathbb{Q}$. It follows from \cite[Proposition 3.19]{lor} that $c(A/\mathbb{Q})=c(E/K)$. Therefore, the abelian variety $A/\mathbb{Q}$ has a $\mathbb{Q}$-rational point of order $p$, has dimension $\frac{p-1}{2}$, and $p \nmid c(A/\mathbb{Q})$. This concludes the proof of our theorem.
\end{proof}

We conclude this section with two examples of Krumm \cite{Krummthesis} which show that there exist elliptic curves $E/K$ defined over cubic (resp. quartic) number fields $K/\mathbb{Q}$ with a $K$-rational point of order $13$ (resp. $13$) and such that $c(E/K)=1$.

\begin{example}\label{example1Krumm} (\cite[Example 5.4.4]{Krummthesis})
Let $K:=\mathbb{Q}[t]/(t^3+2t^2-t-1)$ and let $E/K$ be the elliptic curve given by the Weierstrass equation $$y^2+(-2t^2+2)xy+(-9t^2+2t+4)y=x^3+(-9t^2+2t+4)x^2. $$
Then $K/\mathbb{Q}$ is a cubic number field, $E(K)_{\mathrm{tors}} \cong \mathbb{Z}/13\mathbb{Z}$, and $c(E/K)=1$.
\end{example}

\begin{example}\label{example2Krumm} (\cite[Example 5.5.2]{Krummthesis})
Let $K:=\mathbb{Q}[t]/(t^4-t^3-3t^2+t+1)$ and let $E/K$ be the elliptic curve given by the Weierstrass equation $$y^2+(-6t^3-7t^2+4t+4)xy+(-155t^3-170t^2+109t+74)y=x^3+(-155t^3-170t^2+109t+74)x^2. $$
Then $K/\mathbb{Q}$ is a quartic number field, $E(K)_{\mathrm{tors}} \cong \mathbb{Z}/17\mathbb{Z}$, and $c(E/K)=1$.
\end{example}

\section{Elliptic curves over function fields}\label{sectionfunctionfields}

For this section, we let $k$ be a finite field of characteristic $p>0$. We first recall a few definitions. An elliptic curve $E/k(t)$ is called {\it constant} if there exists an elliptic curve $E_0/k$ such that $E \cong E_0 \times_k k(t)$. An elliptic curve $E/k(t)$ is called {\it isotrivial} if there exists a finite extension $K'/k(t)$ such that the base change $E_{K'}/K'$ of $E/k(t)$ to $K'$ is constant. Finally, an elliptic curve $E/k(t)$ will be called {\it non-isotrivial} if it is not isotrivial.

\begin{proposition}\label{theoremk(t)}
Let $k$ be a finite field of characteristic $p \geq 5$, and let $E/k(t)$ be a non-isotrivial elliptic curve with a $k(t)$-rational point of order $p$. Then $p \mid c(E/k(t))$.
\end{proposition}
\begin{proof} If $k$ is a finite field of characteristic $p$ and there exists an elliptic curve $E/k(t)$ with a point of order $p$, then $p=5,7$, or $11$ (see \cite[Corollary 1.8]{mcdonaldtorsionongenus0}).

Assume that $p=5$ and let $E/k(t)$ be a non-isotrivial elliptic curve with a $k(t)$-rational point of order $5$. The curve $E/k(t)$ can be given by an equation of the form
\begin{align*}
    y^2+(1-f)xy-fy=x^3-fx^2,
\end{align*}
for some non-constant $f \in k(t)$ (see \cite[Table 2]{mcdonaldtorsionongenus0}). The discriminant of this equation is $$\Delta=f^5(f^2-11f-1).$$ Since $f$ is non-constant there exists a valuation $v$ of $k[t]$ such that $v(f)>0$. Using the observation of \ref{tatealgorithm} we see that $E/k(t)$ has split multiplicative reduction modulo $v$ with $5 \mid c_v(E/k(t))$ and, hence, $5 \mid c(E/k(t))$.

Assume that $p=7$ and let $E/k(t)$ be a non-isotrivial elliptic curve with a $k(t)$-rational point of order $7$. The curve $E/k(t)$ can be given by an equation of the form
\begin{align*}
    y^2+(1-a)xy-by=x^3-bx^2,
\end{align*}
with $a=f^2-f$ and $b=f^3-f^2$ for some non-constant $f \in k(t)$ (see \cite[Table 2]{mcdonaldtorsionongenus0}). The discriminant of this equation is $$\Delta=f^7(f-1)^7(f^3-8f^2+5f+1).$$ Since $f$ is non-constant there exists a valuation $v$ of $k[t]$ such that $v(f)>0$. Using \ref{tatealgorithm} we see that $E/k(t)$ has split multiplicative reduction modulo $v$ with $7 \mid c_v(E/k(t))$ and, hence, $7 \mid c(E/k(t))$.

Assume that $p=11$ and let $E/k(t)$ be a non-isotrivial elliptic curve with a $k(t)$-rational point of order $11$. The curve $E/k(t)$ can be given by an equation of the form
\begin{align*}
    y^2+(1-a)xy-by=x^3-bx^2,
\end{align*}
with $$a=\frac{(f+3)(f+5)^2(f+9)^2}{3(f+1)(f+4)^4} \quad \text{and} \quad b=a\frac{(f+1)^2(f+9)}{2(f+4)^3},$$ for some non-constant $f \in k(t)$ (see \cite[Table 14]{mcdonaldtorsionongenus0}). The discriminant of this equation is $$\Delta=\frac{2f^2(f+3)^{11}(f+5)^{11}(f+9)^{11}}{(f+4)^{37}(f+1)}.$$

If $v$ is a valuation of $k(t)$ with $v(f+3)>0$, $v(f+5)>0$, or $v(f+9)>0$, then $v(a)>0$ and $v(b)>0$. Therefore, using \ref{tatealgorithm} we see from the Weierstrass equation of $E/k(t)$ that $E/k(t)$ has split multiplicative reduction at $v$ and moreover $11 \mid v(\Delta)$. This implies that $11 \mid c_v(E/k(t))$. Since $f$ is non constant, there exist valuations $v_1$, $v_2$, and $v_3$ of $k[t]$ such that $v_1(f+3)>0$, $v_2(f+5)>0$, and $v_3(f+9)>0$. Therefore, since $11 \mid c_{v_i}(E/k(t))$ for $i \in \{1,2,3\}$, we find that $11^3 \mid c(E/k(t))$.
\end{proof}

\begin{proposition}\label{splitfunctionfieldtamagawa}
Let $K=\mathbb{F}_q(\mathcal{C})$ be the function field of a smooth, projective, and geometrically irreducible curve $\mathcal{C}/\mathbb{F}_q$, where $q$ is a power of a prime $p \geq 5$. Let $E/K$ be a non-isotrivial elliptic curve with a $K$-rational point of order $p$. If $E/K$ has a place of split multiplicative reduction, then $p \mid c(E/K)$.
\end{proposition}
\begin{proof}
Let $v$ be a place of $K$ such that $E/K$ has split multiplicative reduction modulo $v$. It is enough to show that $p \mid v(\Delta)$, where $\Delta$ is the discriminant of a minimal Weierstrass equation for $E/K$. Since the $j$-invariant $j(E)$ of $E/K$ is equal to $\frac{c_4^3}{\Delta}$ and $v(c_4)=0$, to show that $p \mid c_v(E/K)$, it is enough to show that $p \mid v(j(E))$. The latter follows from the following proposition.
\end{proof}
\begin{proposition}(see \cite[Proposition 7.3]{ulmer2011park})
Let $K=\mathbb{F}_q(\mathcal{C})$ be the function field of a smooth, projective, and geometrically irreducible curve $\mathcal{C}/\mathbb{F}_q$, where $q$ is a power of $p$. Let $E/K$ be a non-isotrivial elliptic curve defined over $K$. Then $E/K$ has a $K$-rational point of order $p$ if and only if $j(E) \in K^p$ and $A(E,\omega)$ is a $(p-1)$-st power in $K^{\times}$, where $A(E,\omega)$ is the Hasse invariant of $E/K$.
\end{proposition}

\begin{proposition}\label{basechange}
Let $K=\mathbb{F}_q(\mathcal{C})$ be the function field of a smooth, projective, and geometrically irreducible curve $\mathcal{C}/\mathbb{F}_q$, where $q$ is a power of $p$. Let $E/K$ be a non-isotrivial elliptic curve with a $K$-rational point of order $p$. Then there exists a finite separable extension $K'/K$ such that $p \mid c(E_{K'}/K')$, where $E_{K'}/K'$ is the base change of $E/K$ to $K'$.
\end{proposition}
\begin{proof}
Since $j(E)$ is non-constant and $K$ is the function field of a smooth, projective, and geometrically irreducible curve, there exists a place $v$ of $K$ such that $v(j(E))<0$ (see \cite[Corollary 1.1.20]{stichtenothbook}). Therefore, the curve $E/K$ has potentially multiplicative reduction at $v$. By the semi-stable reduction theorem for elliptic curves there exists a finite extension of $L/K$ such that the base change $E_L/L$ of $E/K$ to $L$ has semi-stable reduction. After a further finite extension $K'/L$ if necessary (so that the slopes of the tangent lines at the node of the reduced curve are defined over the residue field) we can assume that the base change $E_{K'}/K'$ has split multiplicative reduction modulo a place above $v$. Using Proposition \ref{splitfunctionfieldtamagawa} we find that $p \mid c(E_{K'}/K')$. This proves our proposition.
\end{proof}

\bibliographystyle{plain}
\bibliography{bibliography.bib}

\end{document}